# Topological and metric spaces are full subcategories of the category of simplicial objects of the category of filters.

notes by misha gavrilovich[*]

**Abstract**

We observe that the category of topological space, uniform spaces, and simplicial sets are all, in a natural way, full subcategories of the same larger category, namely the simplicial category of filters; this is, moreover, implicit in the definitions of a topological and uniform space.

We use these embeddings to rewrite the notions of completeness, precompactness, compactness, Cauchy sequence, and equicontinuity in the language of category theory, which we hope might be of use in formalisation of mathematics and tame topology. We formulate some arising open questions.

## 1 Introduction.

In this proposal we define two fully faithful embeddings of the category of topological spaces and that of uniform metric spaces into the category of simplicial objects of the category of filters, and, based on this, use these two functors to reformulate several elementary notions including that of being compact, precompact, complete, a Cauchy sequence, and equicontinuity.

We formulate a number of open questions, largely for the author's own use: this proposal is at a very early stage and it is possible that some of the questions are easy or indeed well-known.

If someone already knows the answers or the relevant literature, the author would be delighted to hear about it.

We hope our reformulations suggest that a number of notions as defined in [Bourbaki, General Topology] may conveniently and concisely be expressed in the language of category theory. It may be worthwhile to express them this way, for two reasons: it may provide a fresh point of view on foundations of topology (tame topology) and it may lead to a development of the language of category theory. It is possible that this may be of use in formalisation of foundations of topology.

---

[*]Comments welcome. `mishap@sdf.org`. I thank Dmitry Krachun, Sergei Ivanov and Vladimir Sosnilo for discussions.

## 2 Main constructions and open questions.

### 2.1 The category of simplicial filters.

We say a topological space is *filtered* iff

($F_I$) any superset of a non-empty open set is open.

We say a subset of a topological space is *big* iff it is non-empty and open.

A *filter* is a filtered topological space $X$ which is not discrete. In a filter the intersection of two big, i.e. non-empty open, subsets is big, i.e. non-empty. Indeed, for any two disjoint open non-empty subsets $U$ and $V$, an arbitrary subset $X$ is the intersection $X = (U \cap V) \cup X = (U \cup X) \cap (V \cup X)$ of two open non-empty subsets $U \cup X$ and $V \cup X$.

Let ꝙ*ilt* be the full subcategory of the category of topological spaces whose objects are filtered topological spaces. Let ꝑ*ilt* be the category with the same objects but maps considered up to being equal almost everywhere, i.e. two continuous maps between filtered topological spaces are considered equal in ꝑ*ilt* iff they coincide on a big subset of the source.

This category has all small limits and colimits and a non-commutative tensor product [Blass,Thm.7]. Limits and colimits are set-wise the same as in *Sets* and the topology is defined as the finest/coarsest filtered topology such that the necessary maps are continuous. Let sꝙ*ilt* be the category of simplicial objects in the category ꝙ*ilt* of filtered topological spaces, i.e. $s\text{ꝙ}ilt = Func(\text{Ord}_{<\omega}^{op}, \text{ꝙ}ilt)$ where $\text{Ord}_{<\omega}$ denotes the category of categories corresponding to finite linear orders

$$\bullet_1 \longrightarrow .. \longrightarrow \bullet_n,\ 0 \le n < \omega.$$

There are two natural functors ꝙ*ilt* ⟶ sꝙ*ilt*:

ɪ : $F \longmapsto (F, F, F, ...)$, identity maps

$E : F \longmapsto (F, F \times F, F \times F \times F, ....)$, face and degeneracy maps are coordinate maps $F^n \longrightarrow F^m, (x_1, ..., x_n) \mapsto (x_{i_1}, ..., x_{i_m})$ where $1 \le i_1 \le i_2 \le ... \le i_m \le n$.

The same considerations apply to ꝑ*ilt* and sꝑ*ilt*.

There are two natural inclusions ɢ : Sets ⟶ ꝙ*ilt*: a set $S$ goes to the filter on $S$ with the unique big subset $S$, and ɪ : Sets ⟶ ꝙ*ilt*: a set $S$ goes to the filter on $S$ where all non-empty subsets are big. This gives three fully faithful embeddings sɢ : sSets ⟶ sꝙ*ilt*, sɢ′ : sSets ⟶ sꝑ*ilt*, and sɪ : sSets ⟶ sꝙ*ilt*.

### 2.2 Topological and metric spaces as simplicial filters

Topological and uniform spaces are defined [Bourbaki, Ch1., Ch.2] as systems of neighbourhood filters satisfying certain compatibility conditions, and lead us to define two fully faithful functors 𝐓 : $Top \longrightarrow s\text{ꝙ}ilt$, 𝓜 : $\mathcal{MU} \longrightarrow s\text{ꝙ}ilt$, and in fact also two fully faithful functors 𝒯 : $Top \longrightarrow s\text{ꝑ}ilt$, 𝓶 : $\mathcal{MU} \longrightarrow s\text{ꝑ}ilt$.

In fact, everything we say below about $\varphi ilt$ and $s\varphi ilt$ holds also for $\vartheta ilt$ and $s\vartheta ilt$, i.e. when maps are considered up to being equal almost everywhere.

In Appendix B we show how to "read off" the latter embedding from the definition of uniform structures in [Bourbaki, Chapter 2].

For a topological space $X$, let $\mathbf{T}(X)$ denote the following object in $s\varphi ilt$.

$$\mathbf{T}: X \longmapsto (|X|, |X| \times |X|, |X| \times |X| \times |X|, ...)$$

with face and degeneracy maps being the coordinate maps

$$|X|^n \longrightarrow |X|^m, (x_1, ..., x_n) \mapsto (x_{i_1}, ..., x_{i_m})$$

where $1 \leq i_1 \leq i_2 \leq ... \leq i_m \leq n$.

A subset $U = \{(x_1, ..., x_n) : (x_1, ..., x_n) \in U\}$ of $|X|^n$ is *big* iff the following formula holds:

$\forall x_1 \in X \exists U_{x_1} \ni x_1$ where $U_{x_1}$ is a neighbourhood of $x_1$
$\forall x_2 \in U_{x_1} \exists U_{x_2} \ni x_2$ where $U_{x_2}$ is a neighbourhood of $x_2$
....
$\forall x_n \in U_{x_{n-1}} \exists U_{x_n} \ni x_n$ where $U_{x_n}$ is a neighbourhood of $x_n$
$(x_1, x_2, ..., x_n) \in U$

Note that any big subset of $|X|^n$ contains the diagonal, in particular the topology on $|X|$ is always antidiscrete. Topology on $X$ is discrete iff the diagonal in $|X|^n, n \geq 2$ is open, equivalently a subset of $|X|^n$ is big iff it contains the diagonal.

For $X$ finite, a subset of $|X|^n, n \geq 2$ is big iff it contains all the non-strictly decreasing sequences in the specialisation preorder, i.e. all the sequences $(x_1, x_2, ..., x_n)$ such that $x_i \in cl_X(x_{i+1}), 1 \leq i < n$.

Note that the topology on $|X|^n$ is the coarsest filter such that the maps $|X|^n \longrightarrow |X|^2, (x_1, ..., x_n) \mapsto (x_i, x_{i+1}), 1 \leq i < n$ are continuous.

For a metric space $M$, let $\mathcal{M}(M)$ denote the following object in $s\varphi ilt$.

$$\mathcal{M}: M \longmapsto (|M|, |M| \times |M|, |M| \times |M| \times |M|, ...)$$

where a subset of $|M|^n$ is *big* iff it contains an $\varepsilon$-neighbourhood of the diagonal $\{(x, .., x) : x \in M\}$. Face and degeneracy maps are coordinate maps

$$|X|^n \longrightarrow |X|^m, (x_1, ..., x_n) \mapsto (x_{i_1}, ..., x_{i_m})$$

where $1 \leq i_1 \leq i_2 \leq ... \leq i_m \leq n$.

Note that, as before, any big subset of $|M|^n$ contains the diagonal, in particular the topology on $|M|$ is always antidiscrete. The diagonal in $|M|^n, n \geq 2$ is open, iff the metric space $M$ is discrete.

Note that the topology on $|M|^n$ is the coarsest filter such that the maps $|X|^n \longrightarrow |X|^2, (x_1, ..., x_n) \mapsto (x_i, x_j), 1 \leq i < j \leq n$ are continuous.

Let $\mathcal{MU}$ denote the category of uniform spaces [Bourbaki, II1.1]. In a similar way $\mathcal{M}(M)$ is defined also for $M$ a uniform space.

Note that permutations of coordinates act on $\mathcal{M}(M)$.

**Remark 1.** The intuition behind these definitions is as follows. [Bourbaki,Introduction] writes: 'a topological structure now enables us to give precise meaning to the phrase "such and such a property holds for all points sufficiently near a": by definition this means that the set of points which have this property is a neighbourhood of a for the topological structure in question.' The notion of a *big subset* enables to concisely express the phrase "such and such a property holds provided a point $a_n$ is sufficiently near $a_{n-1}$ which in turn is sufficiently near $a_{n-1}$, which in turn is sufficently near ..., which is turn is sufficiently near $a_1$". By definition this means that the set of tuples $(a_1, \ldots, a_n)$ of points which have this property is big for the topological structure in question.

Further, [Bourbaki, Introduction] writes: 'As we have already said, a topological structure on a set enables one to give an exact meaning to the phrase "whenever $x$ is sufficiently near $a$, $x$ has the property $P\{x\}$". But, apart from the situation in which a "distance" has been defined, it is not clear what meaning ought to be given to the phrase "every pair of points $x, y$ which are sufficiently near each other has the property $P\{x, y\}$", since a priori we have no means of comparing the neighbourhoods of two different points. Now the notion of a pair of points near to each other arises fre- quently in classical analysis (for example, in propositions which involve uniform continuity). It is therefore important that we should be able to give a precise meaning to this notion in full generality, and we are thus led to define structures which are richer than topological structures, namely *uniform structures*.'

The notion of a *big subset* in a uniform space enables to concisely express the phrase "such and such a property holds provided points $a_1, \ldots, a_n$ are sufficiently near to each other". By definition this means that the set of tuples $(a_1, \ldots, a_n)$ of points which have this property is big for the uniform structure in question.

**Claim 1.** $\mathcal{T}: Top \longrightarrow s\mathcal{F}ilt$, and $\mathcal{U}: \mathcal{MU} \longrightarrow s\mathcal{F}ilt$, and $\mathcal{J}: Top \longrightarrow s\mathcal{P}ilt$, and $\mathcal{U}: \mathcal{MU} \longrightarrow s\mathcal{P}ilt$, are fully faithful functors.

*Proof.* The verification is straightforward and we only consider $\mathcal{T}$. The formula is positive and therefore a superset of a big subset is big. The intersection of two neighbourhoods is a neighbourhood [Bourbaki,I1.2,Ax.$(V_{II})$] and this carries though the quantifiers. Finally, each neighbourhood of a point contains the point, and this implies that big subsets necessarily contain the diagonal and thus form a filter. To see continuity of a degeneracy map $(x_1, ..., x_n) \longmapsto (x_1, .., x_i, x_{i+1}, ..., x_n)$, pick the same neighbourhood twice, $U_{x_i} = U_{x_{i-1}}$ or $U_{x_1} = X$ if $i-1 < 1$; this uses that an open subset is a neighbourhood of each of its points [Bourbaki,I1.2,Ax.$(V_{IV})$] and that $X$ itself is open [Bourbaki,I1.2,Ax.$(V_I)$]. To check continuity of a face map $(x_1, ..., x_n) \longmapsto (x_1, .., x_i, x_i, x_{i+1}, ..., x_n)$, pick $x_{i+1} = x_i$; this uses that a point is contained in any neighbourhood of itself [Bourbaki,I1.2,Ax.$(V_{III})$]. The functor is faithful because the morphism on $X^1$ uniquely determines morphisms on all the other Cartesian powers. A function $f: X \longrightarrow Y$ is continuous iff for each point $x \in X$ and each neighbourhood $V_y$ of $y = f(x)$ there is a neighbourhood $U_x$ such that $V_y \subset f(U_x)$. This is implied by the fact that the preimage of a big set $\{y\} \times V_y \cup \bigcup_{y' \neq y} \{y'\} \times Y$ contains $\{x\} \times U_x$ for some $U_x$, as it is big. $\square$

For a metric space $M$, the functor $\mathcal{N}(M): \mathrm{Ord}_{<\omega}^{op} \longrightarrow s\mathbf{\Phi}ilt$ uniquely factors via $FinSets^{op}$: $\mathrm{Ord}_{<\omega}^{op} \xrightarrow{id} FinSets^{op} \longrightarrow s\mathbf{\Phi}ilt$. And in fact it might be better to view it this way:

**Claim 2.** *To give a uniform space is the same as to give a functor $M : FinSets^{op} \longrightarrow sFilt$ such that*

*$|M(S)| = X^S = Hom_{Sets}(S, X)$, for any finite set $S$ and some fixed set $X$*

*the only non-empty open subset of $M(\{\bullet\})$ is $|M(\{\bullet\}|$ itself*

*the filter on $M(S)$ is the coarsest filter such that the induced maps*

$$M(S) \xrightarrow{M(\{s_i, x_{s+1}\} \longrightarrow S)} M(\{s_i, s_{i+1}\})$$

*are continuous for the inclusion maps $\{s_i, s_{i+1}\} \longrightarrow S$, $1 \leq i \leq \mathrm{card}(S)$, for some, equiv. each, enumeration $S = \{s_1, ..., s_{\mathrm{card}(S)}\}$*

*Proof.* It is straightforward to verify that the construction above gives such a functor. Conversely, such a functor defines a filter on $X \times X$, for some set $X$, and it is straightforward to verify this filter determines a uniform structure on $X$ using the axioms [Bourbaki,II 1.1,Ax.$(U_i) - (U_{III})$]. □

These two embeddings immediately give rise the following questions.

**Question 1.** (Category theory and homotopy)

1. Do functors $\mathcal{T}, \mathcal{N}$ have adjoints? Do they preserve limits and colimits?

2. Characterise the full subcategories $\mathcal{T}(Top)$ and $\mathcal{N}(\mathcal{MU})$ of $s\mathbf{\Phi}ilt$ in terms of the ambient category $s\mathbf{\Phi}ilt$ or $s\mathbf{\Phi}ilt$.

3. Characterise systems of neighbourhoods such that the construction of $\mathcal{T}$ gives rise to a simplicial object.

4. Does a model structure on $Top$ extend to a model structure on $s\mathbf{\Phi}ilt$ or $s\mathbf{\Phi}ilt$? Does $s\mathbf{\Phi}ilt$ or $s\mathbf{\Phi}ilt$ have an interesting model structure?

The following questions are vague.

Everywhere below we may talk of $s\mathbf{\Phi}ilt$ instead of $\mathbf{\Phi}ilt$, and sometimes we omit $s\mathbf{\Phi}ilt$.

**Question 2.** (Naive homotopy theory)

1. What is the "right" notion in $s\mathbf{\Phi}ilt$ or $s\mathbf{\Phi}ilt$ of a real line interval $[0, 1]$, a fibration, and path, loop, and suspension objects?

2. Is there an interesting object of $s\mathbf{\Phi}ilt$ which corresponds to the path space of a topological space and which is more "finitary"? Note that in $s\mathbf{\Phi}ilt$ topological spaces have dimension 2 and that path space is thought of a space "shifted". Can this "shift" be realised in $s\mathbf{\Phi}ilt$ somehow, e.g. so that the $s\mathbf{\Phi}ilt$-path space of a topological space have dimension 3 ?

3. Is there an interesting object of $s\mathbf{\Phi}ilt$ which corresponds to a foliation, particularly an irrational foliation?

**Question 3.** ("Combinatorial" definition of compactness and completeness)

1. A number of elementary topological properties can be defined by, in a sense, combinatorial expressions, by taking iterated orthogonals in $Top$ of a single morphism between finite topological spaces [Gavrilovich, Lifting Property]. Calculate these expressions in $s\mathcal{P}ilt$ using the embedding $\mathcal{T}: Top \longrightarrow s\mathcal{P}ilt$. Note that this would give properties of both topological spaces and metric spaces. Do they define the same properties of topological spaces? Do they provide an interesting analogy between topological spaces and metric spaces, e.g. compactness [Bourbaki, I10.2, Thm.1(d), p.101] and completeness [Bourbaki, II3.6, Prop.11]?

2. The following is an example of a precise conjecture. Evaluated in $Top$, the following expression defines the class of almost(?) all proper maps [Gavrilovich, Lifting Property, Claim 1]:

$$((\{\{o\} \longrightarrow \{o \to c\}\}^r)_{<5})^{lr}$$

If evaluated in $s\mathcal{P}ilt$, does the same expression define both the class of compact space and complete metric space, i.e. is the following true:

   (a) A Hausdorff space $X$ is compact iff
   
   $$\mathcal{T}(X \longrightarrow \{pt\}) \in ((\{\mathcal{T}(\{o\} \longrightarrow \{o \to c\})\}^r)_{<5})^{lr}$$
   
   (b) A metric space $M$ is complete iff
   
   $$\mathcal{T}_l(M \longrightarrow \{pt\}) \in ((\{\mathcal{T}(\{o\} \longrightarrow \{o \to c\})\}^r)_{<5})^{lr}$$

**Question 4.** (Tame topology of Grothendieck)

1. Does this point of view shed light on tame topology of Grothendieck, k, i.e. a foundation of topology "without false problems" and "wild phenomena" "at the very beginning"? Is there a better construction of "the" tubular neighbourhood of closed tame subspace in a tame space? I quote a specific suggestion by Grothendieck [Esquisse d'un Programme, translation,5,p.33]:

   > Among the first theorems one expects in a framework of tame topology as I perceive it, aside from the comparison theorems, are the statements which establish, in a suitable sense, the existence and uniqueness of "the" tubular neighbourhood of closed tame subspace in a tame space (say compact to make things simpler), together with concrete ways of building it (starting for instance from any tame map $X \longrightarrow \mathbb{R}^+$ having $Y$ as its zero set), the description of its "boundary" (although generally it is in no way a manifold with boundary!) $\partial T$, which has in $T$ a neighbourhood which is isomorphic to the product of $T$ with a segment, etc. Granted some suitable equisingularity hypotheses, one expects that $T$ will be endowed, in an essentially unique way, with the structure of a locally trivial fibration over $Y$, with $\partial T$ as a subfibration.

**Question 5.** (History and formalisation of mathematics)

Early works on topology talk about topological spaces in terms of neighbourhood systems. Could it be that they are implicitly trying to express (say, functorial) constructions in $s\mathcal{P}ilt$, are implicitly using category theoretic language but describing it in words? Can this question be made precise?

For example, it was a convention to always mean by $U(x)$ a neighbourhood of a point $x$, and that's somewhat natural from the point of view of our definition of $\mathbf{T}: Top \longrightarrow s\Phi ilt$.

In Appendix B we show how to "read off" our construction of $\mathcal{H}: \mathcal{MU} \longrightarrow s\Phi ilt$ from Bourbaki. Arguably, $\mathbf{T}: Top \longrightarrow s\Phi ilt$ and the combinatorial definitions of elementary topological properties [Gavrilovich, Lifting Property] are "implicitly contained" in Bourbaki. In what sense are these reformulations "contained implicitly" there? Can this sense be made explicit?

Could these reformulations be of use in formalisation of topology and analysis, e.g. Chapter 1 and 2 of Bourbaki?

## 2.3 Elementary theory of topological and metric spaces

We reformulate several notions from [Bourbaki, Chapter 1, 2] in terms of functors $\mathbf{T}: Top \longrightarrow s\Phi ilt$ and $\mathcal{H}: \mathcal{MU} \longrightarrow s\Phi ilt$.

### 2.3.1 Compact and complete metric spaces

An *ultrafilter* is a filter such that if the union of finitely many sets is open, then one of them is; equivalently, each subset $A$ is either open or closed.

With a filter $\mathcal{F}$ on the set of points of a topological space $X$ associate [Bourbaki, II5,Example] a topological space $X \cup_{\mathcal{F}} \{\infty\}$ such that $\mathcal{F}$ is the neighbourhood filter of $\infty$: $|X \cup_{\mathcal{F}} \{\infty\}| = |X| \cup \{\infty\}$, and a subset is open iff it is either an open subset of $X$ or contains $\infty$ and is a union of $\{\infty\}$ and an open subset of $X$ which is also $\mathcal{F}$-open.

Let $X$ be a topological space such that $|X| = |\mathcal{F}|$. An filter $\mathcal{F}$ *converges* on a topological space $X$ iff one of the two equivalent conditions holds [Bourbaki,II7, Def.1]:

there is a point $\infty \in X$ such that each $X$-open neighbourhood of $X$ is also $\mathcal{F}$-open the obvious map $|F| \longrightarrow X$ extends to a map $|F| \cup_{\mathcal{F}} \{\infty\} \longrightarrow X$

**Question 6.** (Ultrafilters and cluster points)

1. Find a category theoretic way to work with ultrafilters, possibly using that in $Top$, for an ultrafilter $\mathcal{F}$, in $Top\ \mathcal{F} \longrightarrow \mathcal{F} \cup_{\mathcal{F}} \{\infty\} \measuredangle g$ for each closed map of finite topological spaces, and, more generally, $\mathcal{F} \longrightarrow \mathcal{F} \cup_{\mathcal{F}} \{\infty\} \in (\{\{0\} \longrightarrow \{0 \to 1\}\}^l_{finite})^r$ (see [Gavrilovich, Lifting property, Claim 1]).

2. Find a category theoretic way or convenient notation to work with cluster points of filters rather than limit point of ultrafilters.

A topological space $X$ is *quasi-compact* iff one of the two equivalent conditions holds [Bourbaki, I10.2, Thm.1(d), p.101]:

each ultrafilter on $X$ converges

for each ultrafilter $\mathcal{F}$ it holds in $Top\ \mathcal{F} \longrightarrow \mathcal{F} \cup_{\mathcal{F}} \{\infty\} \measuredangle X \longrightarrow \{\bullet\}$ for each ultrafilter $\mathcal{F}$ it holds in $s\varphi ilt$ $\mathrm{T}(\mathcal{F}) \longrightarrow \mathrm{T}(\mathcal{F} \cup_{\mathcal{F}} \{\infty\}) \,\measuredangle\, \mathrm{T}(X) \longrightarrow \mathrm{T}(\{\bullet\})$

More generally, a map $X \longrightarrow Y$ is *proper* iff for each ultrafilter $\mathcal{F}$ either of the following equivalent conditions holds [Bourbaki, I10.2, Thm.1(d), p.101]:

in $Top$, $\mathcal{F} \longrightarrow \mathcal{F} \cup_{\mathcal{F}} \{\infty\} \,\measuredangle\, X \longrightarrow Y$

in $s\varphi ilt$, $\mathrm{T}(\mathcal{F}) \longrightarrow \mathrm{T}(\mathcal{F} \cup_{\mathcal{F}} \{\infty\}) \,\measuredangle\, \mathrm{T}(X) \longrightarrow \mathrm{T}(Y)$

A *Cauchy filter* $\mathcal{F}$ on a metric space $M$ ([Bourbaki,II3.1,Def.2]) is a filter on $|M|$ such that one of the two equivalent conditions holds:

for each $\epsilon > 0$ there is a $\mathcal{F}$-open non-empty subset $V \subset |M|$ of diameter at most $\varepsilon$

the map $\mathcal{F} \times \mathcal{F} \longrightarrow |M| \times |M|$ is continuous where $|M| \times |M|$ is equipped with the topology coming from $\mathfrak{N}(M)$

the obvious map $\mathrm{E}(\mathcal{F}) \longrightarrow \mathfrak{N}(M)$ is well-defined

A *Cauchy sequence* in $M$ is a map $\mathrm{E}(\omega_{cofinite}) \longrightarrow \mathfrak{N}(M)$ where $\omega_{cofinite}$ is the set of natural numbers equipped with cofinite topology (i.e. a subset is closed iff it is finite).

A metric space is *precompact* iff one of the two equivalent conditions holds :

for each $\varepsilon > 0$ there is a finite covering of $M$ by subsets of diameter at most $\varepsilon$ [Bourbaki, II4,Thm.3]

each ultrafilter on $M$ is a Cauchy ultrafilter [Bourbaki, II4,Exer.5]

for each ultrafilter it holds in $s\varphi ilt$ $\iota(\mathcal{F}) \longrightarrow \mathrm{E}(\mathcal{F}) \,\measuredangle\, \mathfrak{N}(M) \longrightarrow \mathfrak{N}(\{\bullet\})$

A metric space $M$ is *complete* iff one of the two equivalent conditions holds [Bourbaki,II3.3,Def.3]:

each Cauchy filter on $M$ converges

in $s\varphi ilt$, $\mathrm{E}(\mathcal{F}) \longrightarrow \mathrm{E}(X \cup_{\mathcal{F}} \{\infty\}) \,\measuredangle\, \mathfrak{N}(M) \longrightarrow \mathfrak{N}(\{\bullet\})$

in $s\varphi ilt$, $\mathrm{T}(\mathcal{F}) \longrightarrow \mathrm{T}(X \cup_{\mathcal{F}} \{\infty\}) \,\measuredangle\, \mathfrak{N}(M) \longrightarrow \mathfrak{N}(\{\bullet\})$

**Question 7.** Define the completion of a uniform space [Bourbaki, II3.7] as something like inner hom $\underline{Hom}(\mathrm{E}(\omega_{cofinite}), \mathfrak{N}(M))$. Develop the theory [Bourbaki, II3,4] of complete and precompact uniform spaces in terms of $s\varphi ilt$ and the lifting properties.

Let $M$ be a metric space. The following are equivalent [Bourbaki,II1.2,Def.3]:

topological space $M_{top}$ is homeomorphic to $|M|$ with the topology induced from the metric on $M$

there is an arrow $T(M_{top}) \xrightarrow{\gamma} \mathcal{M}(M)$ and for each topological space $X$, in $s\varphi ilt$ any map $T(X) \longrightarrow \mathcal{M}(M)$ factors as

$$T(X) \longrightarrow T(M_{top}) \xrightarrow{\gamma} \mathcal{M}(M)$$

For a compact space $K$, there exists a unique uniform space $K_{uni}$ which induces on $K$ its topology. In other words, there is a unique map $T(K) \xrightarrow{\gamma} \mathcal{M}(K_{uni})$ such that each map $T(K) \longrightarrow \mathcal{M}(M)$ factors as

$$T(K) \xrightarrow{\gamma} \mathcal{M}(K_{uni}) \longrightarrow \mathcal{M}(M)$$

[Bourbaki, II4.1,Thm.1].

**Remark 2.** The notion of the topology induced by a metric is reminiscent of an adjoint functor to $T$. Does either $\mathcal{M}$ or $T$ have adjoints?

### 2.3.2 Equicontinuous functions and Arzela-Ascoli theorem

Let $X$ be a topological space, let $M$ be a metric space, and let $(f_i)_{i \in \mathbb{N}}$ be a family of functions $f_i : X \longrightarrow M$.

The family $f_i$ is *equicontinuous* if either of the following equivalent conditions holds:

for every $x \in X$ and $\epsilon > 0$, there exists a neighbourhood $U$ of $x$ such that $d_Y(f_i(x'), f_i(x)) \leq \epsilon$ for all $i \in \mathbb{N}$ and $x' \in U$ the map $T(X) \times \mathfrak{t}(\{\mathbb{N}\}) \longrightarrow \mathcal{M}(M)$, $(x, i) \longmapsto f_i(x)$ is well-defined the map $T(X) \times \mathfrak{t}(\mathbb{N}_{cofinite}) \longrightarrow \mathcal{M}(M)$, $(x, i) \longmapsto f_i(x)$ is well-defined

If $X = (X, d_X)$ is also a metric space, we say that the family $f_i$ is *uniformly equicontinuous* iff either of the following equivalent conditions holds:

for every $\epsilon > 0$ there exists a $\delta > 0$ such that $d_Y(f_i(x'), f_i(x)) \leq \epsilon$ for all $i \in \mathbb{N}$ and $x', x \in x$ with $d_X(x, x') \leq \delta$ the map $\mathcal{M}(X) \times \mathfrak{t}(\{\mathbb{N}\}) \longrightarrow \mathcal{M}(M)$, $(x, i) \longmapsto f_i(x)$ is well-defined the map $\mathcal{M}(X) \times \mathfrak{t}(\mathbb{N}_{cofinite}) \longrightarrow \mathcal{M}(M)$, $(x, i) \longmapsto f_i(x)$ is well-defined

The family is *uniformly Cauchy* iff either of the following equivalent conditions holds:

for every $\epsilon > 0$ there exists a $\delta > 0$ and $N > 0$ such that $d_Y(f_i(x'), f_j(x)) \leq \epsilon$ for all $i, j > N$ and $x', x \in x$ with $d_X(x, x') \leq \delta$.

the map $\mathcal{M}(X) \times E(\mathbb{N}_{cofinite}) \longrightarrow \mathcal{M}(M)$, $(x, i) \longmapsto f_i(x)$ is well-defined

Here $\{\mathbb{N}\}$ denotes the trivial filter on $\mathbb{N}$ with a unique big subset $\mathbb{N}$ itself, and $\mathbb{N}_{cofinite}$ denotes the filter of cofinite subsets of $\mathbb{N}$.

**Question 8.** (Arzela-Ascoli)

1. Reformulate various notions of equicontinuity and convergence of a family of functions $f_i : X \longrightarrow M$ in terms of maps in $s\varphi ilt$ using e.g. $\mathbf{k}(\mathbb{N}_{cofinite})$, $\mathrm{E}(\mathbb{N}_{cofinite})$, $\mathbf{k}(\mathbb{N}_{cofinite}\cup_{\mathbb{N}_{cofinite}}\{\infty\})$, $\mathrm{T}(\mathbb{N}_{cofinite}\cup_{\mathbb{N}_{cofinite}}\{\infty\})$, $E(\mathbb{N}_{cofinite}\cup_{\mathbb{N}_{cofinite}}\{\infty\})$, $\mathrm{T}(\mathbb{N}_{cofinite})$, $\mathrm{T}(X)$, $\mathcal{M}(X)$, and $\mathcal{M}(M)$.

2. Reformulate and prove Arzela-Ascoli theorem in terms something like inner $Hom$ in $s\varphi ilt$ and the lifting properties defining precompactness, compactness etc.

3. Define function spaces in terms of something like inner $Hom$ in $s\varphi ilt$.

# 3 Appendix A.

# Embeddings of geometric categories

Here we define several embeddings of geometric categories of metric spaces into the category of "infinitary" simplicial objects of the category of filters, notably the categroy of metric spaces up to quasi-isometry. To make the exposition self-contained, we repeat here some of the notation introduced above. In part this is motivated by a remark in [Gromov, Hyperbolic dynamics, 2.7,p.54, footnote 90].

It may be interesting to consider here $\mathcal{F}ilt$ instead of $\mathcal{F}ilt$.

Let $\mathrm{Ord}_{<\alpha}$ denote the category of finite ordinals less than $\alpha$ and non-decreasing maps; equivalently but more conceptually, this is the full subcategory of the category of categories consisting of the categories $\bullet_0 \longrightarrow \bullet_1 \longrightarrow ...$ corresponding to well-ordered sets of size less than $\alpha$. When $\alpha = \omega + 1$, the category $\mathrm{Ord}_{<\omega}$ is the category of finite ordinals usually denoted $\Delta$.

For a category $C$ and ordinal $\alpha$, $<\alpha$-*simplicial objects in* $C$ is a functor $F: \mathrm{Ord}_{<\alpha}^{op} \longrightarrow C$. These objects naturally form a category which we denote $s_{<\alpha}C = Func(\mathrm{Ord}_{<\alpha}^{op}, C)$ of functors from $\mathrm{Ord}_{<\alpha}^{op}$ to $C$. When $\alpha = \omega + 1$, this is the usual category of simplicial objects of $C$.

With an object $X$ we can associate two $<\alpha$-simplicial objects in $C$ as follows. $\iota(X)$ sends each ordinal to $X$ itself and each morphism to the identity The functor $E_\alpha(X)$ sends an ordinal $\beta < \alpha$ to the Cartesian power $X^\beta$, and morphisms are sent to the coordinate maps.

These two functors define two fully faithful embeddings of $C$ into $s_{<\alpha}C$. $\iota: C \longrightarrow s_{<\alpha}C$ and $E: C \longrightarrow s_{<\alpha}C$.

Let $\mathcal{F}ilt$ be the category of filters, i.e. the full subcategory of the category of topological spaces consisting of spaces such that any superset of a non-empty open set is open.

## 3.1 Metric spaces as "infinitary" simplicial filters

We define several embeddings of categories of metric spaces with various kinds of geometric maps, e.g. uniformly continuous maps, Lipschitz maps on large scale. We do so by definition various filters on (possibly infinite) Cartesian powers of a metric space which preserve certain geometric information about the metric space. In the usual way these collections of filters give rise to simplicial objects of $s_{\leq \omega}\mathcal{F}ilt$.

Let $M$ be a metric space. Let us now define a number of topologies on Cartesian powers of $|M|$.

> A non-empty subset of $M^n$ is $\tau$-open (big) iff the following formula holds:
> $\forall x_1 \in M \exists U_{x_1} \ni x_1$ where $U_{x_1}$ is a neighbourhood of $x_1$
> $\forall x_2 \in U_{x_1} \exists U_{x_2} \ni x_2$ where $U_{x_2}$ is a neighbourhood of $x_2$
> ....

$\forall x_n \in U_{x_{n-1}} \exists U_{x_n} \ni x_n$ where $U_{x_n}$ is a neighbourhood of $x_n$
$(x_1, x_2, ..., x_n) \in U$

A non-empty subset of $M^n$ is $\tau_U$-open iff it contains an $\varepsilon$-neighbourhood of the diagonal $\{(x, x, ..., x) : x \in |M|\}$ for some $\varepsilon > 0$, i.e. $U \subseteq |M|^n$ is open iff there is $\varepsilon > 0$ such that for each $x_1, ..., x_n \in M$, it holds $(x_1, ..., x_n) \in U$ provided there is $x \in M$ such that $dist(x, x_i) < \varepsilon$, $i = 1, .., n$.

Note that no proper subset $\emptyset \subsetneq U \subsetneq |M|$ is open as you may take $x = x_1$.

Fix a real number $D > 0$. A non-empty subset of $M^n$ is $\tau_D$-open iff it contains an $D$-neighbourhood of the diagonal $\{(x, x, ..., x) : x \in |M|\}$ for some $\varepsilon > 0$, i.e. $U \subseteq |M|^n$ is open iff there for each $x_1, ..., x_n \in M$, it holds $(x_1, ..., x_n) \in U$ provided there is $x \in M$ such that $dist(x, x_i) \leq D$, $i = 1, .., n$.

Note that no proper subset $\emptyset \subsetneq U \subsetneq |M|$ is open as you may take $x = x_1$.

A non-empty subset of $M^\omega$ is $\tau_\mathcal{L}$-open iff there is $\lambda > 0$, $N > 0$, $D < \lambda N$ such that for each $x_1, ..., x_n, ... \in M$, it holds $(x_1, ..., x_n, ...) \in U$ provided there is $x \in M$ such that $x = x_1 = ... = x_N$ and $dist(x, x_i) \leq \lambda i - D$ for each $i > N$.

A non-empty subset $U$ of $M^n$ is $\tau_\mathcal{L}$-open iff it contains the diagonal $\{(x, .., x) : x \in M\}$.

A non-empty subset of $M^\omega$ is $\tau_{\mathcal{L}_1}$-open iff there is $N > 0$, $D < N$ such that for each $x_1, ..., x_n, ... \in M$, it holds $(x_1, ..., x_n, ...) \in U$ provided there is $x \in M$ such that $x = x_1 = ... = x_N$ and $dist(x, x_i) \leq i - D$ for each $i > N$.

A non-empty subset $U$ of $M^n$ is $\tau_1$-open iff it contains the diagonal $\{(x, .., x) : x \in M\}$.

A map $f : |M| \longrightarrow |N|$ induces a map $f_n : |M|^n \longrightarrow |N|^n$. The following is easy to check:

For $n > 1$, $f_n$ is is $\tau$-continuous iff it is continuous.

For $n > 1$, $f_n$ is is $\tau_U$-continuous iff it is uniformly continuous.

For $n > 1$, $f_n$ is $\tau_D$-continuous iff for each $x, y \in M$ $dist(x, y) \leq D$ implies $dist(f(x), f(y)) \leq D$

$f_\omega$ is $\tau_\mathcal{L}$-continuous iff it is $\lambda$-Lipschitz on large scale for some $\lambda, D > 0$, i.e. $dist(f(x), f(y)) \leq \lambda dist(x, y)$ whenever $dist(x, y) \geq D$, $x, y \in M$.

$f_\omega$ is $\tau_{\mathcal{L}_1}$-continuous iff it is 1-Lipschitz on large scale, i.e. for some $D$ for each $x, y \in M$ $dist(f(x), f(y)) \leq dist(x, y) + D$

A map $f : M \longrightarrow M$ is *an almost isometry* iff either of the following equivalent conditions holds:

$dist(f(x), f(y)) \leq dist(x, y) + D$ for some $D$

$f_\omega : M^\omega \longrightarrow M^\omega$ is $\tau_{\mathcal{L}_1}$-continuous

A map $f : M \longrightarrow M$ is *a quasi-isometry* iff either of the following equivalent conditions holds:

$dist(f(x), f(y)) \leq \lambda dist(x,y) + D$ for some $\lambda, D > 0$ for each $x, y \in M$

$f_\omega : M^\omega \longrightarrow M^\omega$ is $\tau_{\mathcal{L}}$-continuous

A verification shows that these topologies define fully faithful functors

$$m_\mathcal{U} : \mathcal{MU} \longrightarrow Func(\mathrm{Ord}^{\mathrm{op}}_{<\omega}, \textit{Pilt}),$$

$$m_D : \mathcal{M}_D \longrightarrow Func(\mathrm{Ord}^{\mathrm{op}}_{<\omega}, \textit{Pilt}),$$

$$m_\mathcal{L} : \mathcal{ML} \longrightarrow Func(\mathrm{Ord}^{op}_{<\omega+1}, \textit{Pilt}),$$

$$m_{\mathcal{L}_1} : \mathcal{ML}_1 \longrightarrow Func(\mathrm{Ord}^{op}_{<\omega+1}, \textit{Pilt})$$

from the relevant geometric categories of metric spaces.

# 4 Appendix B.

# Reading Bourbaki definition of the uniform spaces

(Bourbaki, II1.1.1) treats metric spaces as uniform spaces; we observe that the uniform space is a simplicial object.

We quote (Bourbaki, I6.1.1) and (Bourbaki, II1.1.1):

> DEFINITION I. A filter on a set $X$ is a set $\mathcal{F}$ of subsets of $X$ which has the following properties:
> $(F_I)$ Every subset of $X$ which contains a set of F belongs to $\mathcal{F}$.
> $(F_{II})$ Every finite intersection of sets of $\mathcal{F}$ belongs to $\mathcal{F}$.
> $(F_{III})$ The empty set is not in $\mathcal{F}$.
>
> 1. DEFINITION OF A UNIFORM STRUCTURE
>
> DEFINITION I. *A uniform structure (or uniformity) on a set $X$ is a structure given by a set $\mathfrak{U}$ of subsets of $X \times X$ which satisfies axioms $(F_I)$ and $(F_{II})$ of Chapter I, 6, no. I and also satisfies the following axioms:*
> $(U_I)$ *Every set belonging to $\mathfrak{U}$ contains the diagonal $\Delta$.*
> $(U_{II})$ *If $V \in \mathfrak{U}$ then $V^{-1} \in \mathfrak{U}$.*
> $(U_{III})$ *For each $V \in \mathfrak{U}$ there exists $W \in \mathfrak{U}$ such that $W \circ W \subset V$.*
> *The sets of $\mathcal{U}$ are called entourages of the uniformity defined on $X$ by $\mathcal{U}$. A set endowed with a uniformity is called* a uniform space. *If $V$ is an entourage of a uniformity on $X$, we may express the relation $(x, x') \in V$ by saying that "$x$ and $x'$ are $V$-close".*

The set of points of a metric space $X$ carries a canonical uniform space: $V \in \mathfrak{U}$ iff $\{(x, x') : dist(x, x') < \varepsilon\} \subset V$ for some $\varepsilon > 0$.

Let us translate the definitions above to the language of arrows: we shall see that a uniform space may be viewed as a simplicial object of the category of topological spaces.

First notice that a filter can equivalently be defined as a non-discrete topology such that a superset of a non-empty open set is necessarily open: a filter $\mathcal{F}$ on a set $X$ defines a topology on $X$ where a subset is open iff it is either $\mathcal{F}$-big or empty. Indeed, Axioms $(F_I)$ and $(F_{II})$ of a filter imply that the family of subsets $\mathfrak{U} \cup \{\varnothing\}$ is a topology on a set $X$.

In this way a uniform structure on a set $X$ defines a topology on $X \times X$.

Axiom $(U_I)$ implies that the diagonal map $X \xrightarrow{(x,x)} X \times X$ is continuous as a map from the set $X$ equipped with antidiscrete topology to the set $X \times X$ equipped with the topology above, and is almost equivalent to this. Indeed, the latter says that an $\mathfrak{U}$-big subset of $X \times X$ either contains the diagonal or does not intersect it.

Axiom $(U_{II})$ says that permuting the coordinates $X \times X \longrightarrow X \times X, (x_1, x_2) \mapsto (x_2, x_1)$ is continuous in this topology.

Define topology on the set $X \times X \times X$ via the pullback square in the category of filter topological spaces $s\mathbf{\Phi}ilt$:

$$\begin{array}{ccc} X \times X \times X \stackrel{(p_1 \times p_2)}{\underset{(p_2 \times p_3)}{\longrightarrow}} & \stackrel{(p_1 \times p_2)}{\longrightarrow} & X \times X^{p_2} \\ \downarrow & & \downarrow \\ X \times X_{p_2} & \longrightarrow & X \end{array}$$

Axiom $(U_{III})$ says that the map $X \times X \times X \xrightarrow{(p_1, p_3)} X \times X, (x_1, x_2, x_3) \mapsto (x_1, x_3)$ is continuous in this topology. Indeed, by definition

$$W_1 \circ W_2 = \{(x_1, x_3) : (x_1, x_2) \in W_1, \ (x_2, x_3) \in W_2\}$$

and the sets of form

$$\{(x_1, x_2, x_3) : (x_1, x_2) \in W_1, \ (x_2, x_3) \in W_2\} = W_1 \times X \cap X \times W_2,$$

$W_1, W_2 \in \mathfrak{U}$, form a base of the pullback topology on $X \times X \times X$. Hence, $(U_{III})$ says that the preimage of an open subset of $X \times X$ under $(p_1, p_3)$ contains an open subset of $X \times X \times X$, i.e. is open (as pullback is taken among filter topologies).

Axiom $(U_I)$ implies that the diagonal map $X \xrightarrow{(x,x)} X \times X$ is continuous as a map from the set $X$ equipped with antidiscrete topology to the set $X \times X$ equipped with the topology above.

Note that $W \circ W$ intersects the diagonal and the continuity of the diagonal map $X \xrightarrow{(x,x)} X \times X$ implies $W \circ W$ contains the diagonal. Thus, in presence of $(U_{III})$, $(U_I)$ is equivalent to the continuity of the diagonal map $X \xrightarrow{(x,x)} X \times X$ in the topologies indicated.

Let $X_1$ denote the set $X$ equipped with the antidiscrete topology. Let $X_2$ and $X_3$ denote the sets $X \times X$ and $X \times X \times X$ equipped with the topologies above. For $n > 3$, let $X_n$ be the pullback in $s\mathbf{\Phi}ilt$

$$\begin{array}{ccc} X_n \stackrel{(p_2 \times \ldots \times p_n)}{\underset{(p_1 \times p_2)}{\longrightarrow}} & \stackrel{(p_2 \times \ldots \times p_n)}{\longrightarrow} & X_{n-1}{}^{p_2} \\ \downarrow & & \downarrow \\ X_2{}_{p_2} & \longrightarrow & X \end{array}$$

The axioms above ensure that the "set-theoretic" face and degeneracy maps

$$(p_{i_1}, \ldots, p_{i_k}) : X \times \ldots \times X \longrightarrow X \times \ldots \times X$$

are continuous. Thus we see that a uniform structure on a set $X$ defines a simplicial complex $X_n$ in $s\mathbf{\Phi}ilt$,

$$(p_{i_1}, \ldots, p_{i_k}) : X_n \longrightarrow X_m$$

**Claim 3.** *A uniform structure on a set $X$ is a simplicial object $X_\cdot$ in the subcategory $s\varphi ilt$ of filter topological spaces equipped with an involution $i : X_\cdot \longrightarrow X_\cdot$ such that*

$X_1$ *is the set $X$ equipped with antidiscrete topology*

*the underlying set of $X_2$ is $X \times X$*

$i : X_\cdot \longrightarrow X_\cdot$ *is the involution permuting the coordinates on $X \times X$*

*for $n > 2$, $X_n$ is the pullback as described above*

**Question 9.** Find a categorical description of the simplicial objects obtained from uniform spaces.

Acknowledgements. To be written.

This work is a continuation of [DMG]; early history is given there. I thank M.Bays, D.Krachun, K.Pimenov, V.Sosnilo, S.Synchuk and P.Zusmanovich for discussions and proofreading; I thank L.Beklemishev, N.Durov, S.V.Ivanov, S.Podkorytov, A.L.Smirnov for discussions. I also thank several students for encouraging and helpful discussions. Chebyshev laboratory, St.Petersburg State University, provided a coffee machine and an excellent company around it to chat about mathematics. Special thanks are to Martin Bays for many corrections and helpful discussions. Several observations in this paper are due to Martin Bays. I thank S.V.Ivanov for several encouraging and useful discussions; in particular, he suggested to look at the Lebesque's number lemma and the Arzela-Ascoli theorem. A discussion with Sergei Kryzhevich motivated the group theory examples.

Much of this paper was done in St.Petersburg; it wouldn't have been possible without support of family and friends who created an excellent social environment and who occasionally accepted an invitation for a walk or a coffee or extended an invitation; alas, I made such a poor use of it all.

This note is elementary, and it was embarrassing and boring, and embarrassingly boring, to think or talk about matters so trivial, but luckily I had no obligations for a time.